\RequirePackage{ifpdf}
\ifpdf 
\documentclass[pdftex]{sigma}
\else
\documentclass{sigma}
\fi

\newtheorem{Theorem}{Theorem}[section]
\newtheorem{Lemma}[Theorem]{Lemma}
\newtheorem{Proposition}[Theorem]{Proposition}
\newtheorem{Corollary}[Theorem]{Corollary}
{\theoremstyle{definition}
\newtheorem{Remark}[Theorem]{Remark}
\newtheorem{Definition}[Theorem]{Definition}
\newtheorem{Example}[Theorem]{Example}
}

\numberwithin{equation}{section}

\usepackage[all]{xy}

\newcommand{\deff}[1]{\textbf{\emph{#1}}}
\newcommand{\func}[3]{#1 \colon #2 \to #3}
\newcommand{\hookfunc}[3]{#1 \colon #2 \hookrightarrow #3}
\newcommand{\eval}[2]{\left. #1 \right|_{#2}}
\newcommand{\pont}[1]{\mathbb{T} {#1}}
\newcommand{\cpont}[1]{\mathbb{T}_{\mathbb{C}} {#1}}
\newcommand{\tang}{T}
\newcommand{\iprod}[2]{\left\langle {#1} , {#2} \right\rangle}
\newcommand{\iiprod}[2]{\left\langle \! \left\langle {#1} , {#2} \right\rangle \! \right\rangle}
\newcommand{\litlet}{\renewcommand{\labelenumi}{(\alph{enumi})}}

\begin{document}

\allowdisplaybreaks

\renewcommand{\PaperNumber}{081}

\FirstPageHeading

\ShortArticleName{Singular Reduction of Generalized Complex Manifolds}

\ArticleName{Singular Reduction of Generalized Complex Manifolds}

\Author{Timothy E.~GOLDBERG}

\AuthorNameForHeading{T.E.~Goldberg}

\Address{Donald and Helen Schort School of Mathematics and Computing Sciences,\\ Lenoir--Rhyne University, Hickory, North Carolina 28601, USA}
\Email{\href{mailto:timothy.goldberg@lr.edu}{timothy.goldberg@lr.edu}}
\URLaddress{\url{http://mat.lr.edu/faculty/goldberg}}

\ArticleDates{Received March 24, 2010, in f\/inal form October 06, 2010;  Published online October 09, 2010}

\Abstract{In this paper, we develop results in the direction of an analogue of Sjamaar and Lerman's singular reduction of Hamiltonian symplectic manifolds in the context of reduction of Hamiltonian generalized complex manifolds (in the sense of Lin and Tolman).  Specif\/ically, we prove that if a compact Lie group acts on a generalized complex manifold in a Hamiltonian fashion, then the partition of the global quotient by orbit types induces a~partition of the Lin--Tolman quotient into generalized complex manifolds.  This result holds also for reduction of Hamiltonian generalized K\"ahler manifolds.}


\Classification{53D20; 53D18; 53C15}

\section{Introduction}

Generalized complex geometry was introduced by Hitchin in \cite{hitchin}, and further developed by his student Gualtieri in his doctoral thesis~\cite{gualtieri}.  It serves as a common ground in which symplectic, Poisson, and complex geometry can meet.  For this reason, there has been much ef\/fort to import ideas and techniques from these other f\/ields into the generalized complex setting.  In particular, many constructions and results from equivariant symplectic geometry have found useful analogues here.  One example is that of Hamiltonian group actions and moment maps, developed in \cite{lintolman}.  (Similar constructions were developed and examined by other groups, such as \cite{MR2323543, MR2378459} and \cite{MR2534281}, but this paper expands specif\/ically on the work of Lin and Tolman.)

Lin and Tolman's construction generalizes the usual symplectic def\/inition, and they go on to prove that one can reduce a generalized complex manifold by its Hamiltonian symmetries, in perfect parallel to Marsden--Weinstein symplectic reduction~\cite{reduction}  (sometimes also credited to Meyer~\cite{meyer}).  Just as in the symplectic case, in order to ensure that the generalized reduced space is a manifold, one must make an assumption regarding freeness of the group action.

In \cite{stratified}, Lerman and Sjamaar proved that if the symplectic reduced space is not a manifold, then the orbit type stratif\/ication of the original symplectic manifold induces the structure of a~\deff{stratified space} (see Def\/inition~1.7 of that paper) on the reduced space, each stratum of which is naturally a symplectic manifold.  The main result of this paper, Theorem~\ref{MainThm}, is a~f\/irst step in the direction of an analagous result for the case of a~singular generalized complex reduced space.  It states that the singular generalized complex reduced space can be partitioned into disjoint generalized complex manifolds.  It is not yet known whether the reduced space in this situation is actually a stratif\/ied space.  A~similar, although distinct, situation was studied in~\cite{SingularDiracReduction}, in which the authors considered the singular reduction of Dirac manifolds.  They analyzed the global quotient of a Dirac manifold by a proper group action as a dif\/ferential space, as in~\cite{MR1848504}, and obtained conditions that guarantee the Dirac structure will descend to the quotient space.

An interesting dif\/ference between the symplectic and generalized complex situations is that the generic result of symplectic reduction is a space with at worst orbifold singularities, whereas for the reduction of a \emph{twisted} generalized complex manifold, the generic result may be a space with worse-than-orbifold singularities.  See Remark~\ref{rem3} below.

Section~\ref{section1} is a rapid introduction to some essential notions from generalized complex geo\-met\-ry.  Section~\ref{section2} reviews some important facts about equivariant cohomology and the orbit type stratif\/ication of $G$-spaces.  Section~\ref{section3} consists of a summary of Hamiltonian actions and reduction in generalized complex geometry.  Finally, Section~\ref{section4} contains the full statement and proof of our main theorem.

An earlier version of this work appeared in the author's doctoral thesis~\cite{MyThesis}, where many def\/initions and calculations are explained in great detail.

Throughout, we use the abbreviations ``GC'' for ``generalized complex'' and ``GK'' for ``gene\-ralized K\"ahler''.  Also, we typically use the same notation and nomenclature to refer to both a~map and its complex linear extension.  Finally, we make use of the \emph{musical notation} for the map between a vector space and its dual induced by a bilinear form.  If $\func{B}{V \times V}{\mathbb{R}}$ is a bilinear form on a real vector space $V$, then we will denote by $\func{B^{\flat}}{V}{V^*}$ the map $v \mapsto \iota_v B := B(v, \cdot)$, where $\iota_v$ denotes the interior product by $v$.  If $B$ is non-degenerate, then $B^{\flat}$ is invertible and we denote its inverse by $B^{\sharp} := \left(B^{\flat}\right)^{-1}$.  We also use the musical notation for vector bundles, sections of their second symmetric powers, and the associated bundle maps.

\section{Generalized complex geometry}\label{section1}

We begin by giving several standard def\/initions and results from generalized complex geometry, which can be found in \cite{gualtieri} or \cite{benbassat}.

For any smooth manifold $M$, the \deff{Pontryagin bundle}, or \deff{generalized tangent bundle}, of $M$ is $\pont{M} := \tang M \oplus \tang^* M$.  This vector bundle carries a natural non-degenerate symmetric metric $\iiprod{\cdot}{\cdot}$ of signature $(n,n)$, def\/ined by
\[
  \iiprod{u+\alpha}{v+\beta}
  :=
  \frac{1}{2} \left( \alpha(v) + \beta(u) \right)
\]
for all $x \in M$ and $u+\alpha, v+\beta \in \mathbb{T}_x M$.  We will use the same notation for the complex bilinear extension of this metric to the complexif\/ication $\cpont{M} := \pont{M} \otimes_{\mathbb{R}} \mathbb{C}$.
These metrics will henceforth be referred to as the \deff{standard metrics} on $\pont{M}$ and $\cpont{M}$.

\begin{Proposition} \label{prop4}
Let $M$ be a manifold.
There is a natural bijective correspondence between the following two structures.
  \begin{enumerate}\itemsep=0pt
    \item Complex linear subbundles $E \subset \cpont{M}$ over $M$ such that $E \cap \overline{E} = 0$ and $E$ is maximally isotropic with respect to the standard metric on $\cpont{V}$.  $($Here $0$ denotes the image of the zero section of $\cpont{M} \to M$, as is customary.$)$
    \item Bundle automorphisms $\mathcal{J}$ of $\pont{M}$ over the identity $M \to M$ such that $\mathcal{J}^2 = - \mathrm{id}$ and $\mathcal{J}$ is orthogonal with respect to the standard metric on $\pont{M}$.
  \end{enumerate}
\end{Proposition}

\begin{Definition}
Let $M$ be a manifold.  Either of the equivalent structures described in Proposition~\ref{prop4} will be called an \deff{almost GC structure} on $M$.  If $M$ is equipped with an almost GC structure $\mathcal{J}$, then $(M,\mathcal{J})$ is an \deff{almost GC manifold}.

Let $E \subset \cpont{M}$ be an almost GC structure on $M$, and for each $x \in M$ let $\func{\pi_x}{\mathbb{T}_{\mathbb{C},x} M}{\tang_{\mathbb{C},x} M}$ be the projection.  The \deff{type} of this almost GC structure at the point $x \in M$ is the complex codimension of $\pi_x(E_x)$ in $\tang_{\mathbb{C},x} M$:
\[
  \operatorname{type}(E)_x
  =
  \dim_{\mathbb{C}} \tang_{\mathbb{C},x} M - \dim_{\mathbb{C}} \pi_x(E_x).
\]

Let $\mathcal{J}_1$ and $\mathcal{J}_2$ be commuting almost GC structures on $M$.  Then $G := -\mathcal{J}_1 \circ \mathcal{J}_2$ is an orthogonal and involutive bundle map $\pont{M} \to \pont{M}$, and there is an associated bilinear form def\/ined by
\[
  (\mathcal{X},\mathcal{Y}) \mapsto \iiprod{G(\mathcal{X})}{\mathcal{Y}}
\]
for all $\mathcal{X}, \mathcal{Y} \in \pont{M}$ in a common f\/iber.
We call $G$ \deff{positive definite} if its associated bilinear form is positive def\/inite, i.e. if $\iiprod{G(\mathcal{X})}{\mathcal{X}} > 0$ for all nonzero $\mathcal{X} \in \pont{M}$.
An  \deff{almost GK structure} on $M$ is a pair of commuting almost GC structures $(\mathcal{J}_1,\mathcal{J}_2)$ on $M$ such that $G := -\mathcal{J}_1 \circ \mathcal{J}_2$ is positive def\/inite.
\end{Definition}

\begin{Remark}
Let $M$ be a manifold.  A maximally isotropic linear subbundle of $\pont{M}$, respectively~$\cpont{M}$, is called a \deff{Dirac structure}, respectively \deff{complex Dirac structure} on $M$.
Thus an almost GC structure on $M$ is a complex Dirac structure $E \subset \cpont{M}$ satisfying $E \cap \overline{E} = \{0\}$.
\end{Remark}

\begin{Definition}
Let $M$ be a manifold, and let $B \in \Omega^2(M)$, where $\Omega^2(M)$ denotes the space of dif\/ferential two-forms on $M$.  The \deff{$B$-transform} of $\pont{M}$ def\/ined by $B$ is the map
\[
  \func{\mathrm{e}^{B}}{\pont{M}}{\pont{M}},
  \qquad
  \mathrm{e}^B := \begin{pmatrix} 1 & 0 \\ B^{\flat} & 1 \end{pmatrix}.
\]
The $B$-transform $\mathrm{e}^B$ is called \deff{closed} or \deff{exact} if the two-form $B$ is closed or exact, respectively.
\end{Definition}

\begin{Proposition}
Let $M$ be a manifold and let $B \in \Omega^2(M)$.
The $B$-field transform $\mathrm{e}^B$ is orthogonal with respect to the standard metrics on $\pont{M}$ and $\cpont{M}$.
It transforms almost GC structures on $V$ by
\[
  \mathcal{J} \mapsto \mathrm{e}^B \circ \mathcal{J} \circ \mathrm{e}^{-B}
  \qquad \text{and} \qquad
   E \mapsto \mathrm{e}^B(E)
\]
for an almost GC structure given equivalently by a map $\mathcal{J}$ or a Dirac structure $E$, and it preserves types.  It transforms almost GK structures $(\mathcal{J}_1,\mathcal{J}_2)$ by transforming $\mathcal{J}_1$ and $\mathcal{J}_2$ individually.
\end{Proposition}

The Lie bracket def\/ines a skew-symmetric bilinear bracket on sections of the tangent bund\-le~$\tang M$.  This can be extended to a skew-symmetric bilinear bracket on sections of the Pontryagin bundle $\pont{M}$, called the \deff{Courant bracket}, def\/ined by
  \[
    [ X+\alpha, Y+\beta ]
    :=
    [X,Y] + \mathcal{L}_X \beta - \mathcal{L}_Y \alpha
    - \frac{1}{2} \, \mathrm{d} \left( \beta(X) - \alpha(Y) \right)
  \]
for all $X+\alpha, Y+\beta \in \Gamma(\pont{M})$, where $\Gamma(\pont{M})$ denotes the space of smooth sections of $\pont{M} \to M$.
Here the bracket on the right-hand side is the usual Lie bracket of vector f\/ields, and $\mathcal{L}$ denotes Lie dif\/ferentiation. For each closed dif\/ferential three-form $H \in \Omega_{\text{cl}}^3(M)$, there is also the \deff{$H$-twisted Courant bracket}, def\/ined by
  \[
    [ X+\alpha, Y+\beta ]_H
    :=
    [ X+\alpha, Y+\beta ]
    + \iota_Y \iota_X H
  \]
$X+\alpha, Y+\beta \in \Gamma(\pont{M})$.
Both the Courant and the $H$-twisted Courant brackets extend complex linearly to brackets on smooth sections of the complexif\/ied Pontryagin bundle $\cpont{M}$, which will be denoted the same way.

\begin{Definition}
Let $M$ be a manifold, and let $L$ be a real (respectively complex) linear subbundle of $\pont{M}$ (respectively $\cpont{M}$).  Then $L$ is \deff{Courant involutive} if the space $\Gamma(L)$ of smooth sections of $L$ is closed under the Courant bracket, i.e. $\big[ \Gamma(L), \Gamma(L) \big] \subset \Gamma(L)$.  If $H \in \Omega_{\text{cl}}^3(M)$, we similarly def\/ine \deff{$H$-twisted Courant involutive}.

Let $E \subset \cpont{M}$ be an almost GC structure on $M$.  This is a \deff{GC structure} if $E$ is Courant involutive, in which case $(M,E)$ is a \deff{GC manifold}.  If $H \in \Omega_{\text{cl}}^3(M)$ and $E$ is $H$-twisted Courant involutive, then $E$ is an \deff{$H$-twisted GC structure}, and $(M,E,H)$ is a \deff{twisted GC manifold}.

Let $(\mathcal{J}_1, \mathcal{J}_2)$ be an almost GK structure on $M$.  This is a \deff{GK structure} if both $\mathcal{J}_1$ and~$\mathcal{J}_2$ are Courant involutive, in which case $(M,\mathcal{J}_1,\mathcal{J}_2)$ is a \deff{GK manifold}.  If $H \in \Omega_{\text{cl}}^3(M)$ and~$\mathcal{J}_1$ and $\mathcal{J}_2$ are $H$-twisted Courant involutive, then this is an \deff{$H$-twisted GK structure}, and $(M,\mathcal{J}_1,\mathcal{J}_2,H)$ is a \deff{twisted GK manifold}.
\end{Definition}

\begin{Remark}
Let $M$ be a manifold and $D$ be a real or complex Dirac structure on $M$.  Then $D$ is called \deff{closed}, or \deff{integrable}, if the space $\Gamma(D)$ of smooth sections of $D$ is Courant involutive.  Thus a GC structure on $M$ is a closed complex Dirac structure $E \subset \cpont{M}$ such that $E \cap \overline{E} = 0$.
\end{Remark}

\begin{Proposition}[Proposition 3.23 of \cite{gualtieri}]
Let $M$ be a manifold, let $H \in \Omega_{\text{cl}}^3(M)$, and let $B \in \Omega^2(M)$.  The $B$-transform of an $H$-twisted GC structure is an ($H + \mathrm{d} B$)-twisted GC structure.  Thus, a closed $B$-transform of an untwisted GC structure is untwisted.
\end{Proposition}

\begin{Example} \label{example4} \qquad
\begin{enumerate}\itemsep=0pt
\item Let $(M,\omega)$ be an \deff{almost symplectic manifold}, meaning that $\omega \in \Omega^2(M)$ is a non-degenerate form on $M$, but not necessarily closed.  This def\/ines an almost GC structure~$\mathcal{J}_{\omega}$ on~$M$ by
\[
  \mathcal{J}_{\omega} :=
  \begin{pmatrix}
    0 & - \omega^{\sharp} \\
    \omega^{\flat} & 0
  \end{pmatrix}
\] of type $0$ at every point.
It has associated Dirac structure def\/ined by
\[
  E_{\omega, x}
  =
  \{ X - i \, \omega^{\flat}_x(X)
  \mid X \in \tang_{\mathbb{C}, x} M \},
\]
for each $x \in M$.  As discussed in Section 3 of \cite{gualtieri}, it is a~GC structure if and only if $\mathrm{d} \omega = 0$, i.e.\ if and only if~$\omega$ is a symplectic structure on~$M$.

\item Let $(M,I)$ be an \deff{almost complex manifold}, meaning that $I^2 = -\mathrm{id}_{\tang M}$ but $I$ is not necessarily integrable.  This def\/ines an almost GC structure $\mathcal{J}_{I}$ on $M$ by
\[
  \mathcal{J}_{I}
  :=
  \begin{pmatrix}
  -I & 0 \\
  0 & I^*
  \end{pmatrix}.
\]
of type $n$ at every point.  It has associated Dirac structure def\/ined by
\[
  E_{I}
  = \tang_{0,1} M \oplus \tang_{1,0}^* M,
\]
where $\tang_{1,0} M, \tang_{0,1} M \subset \tang_{\mathbb{C}} M$ denote the $\pm i$-eigenbundles of $I$.
 As discussed in Section 3 of \cite{gualtieri}, it is a GC structure if and only if $I$ is integrable, i.e. if and only if $I$ is a complex structure on $M$.

\item Let $M$ be a \deff{K\"ahler manifold} with K\"ahler form $\omega \in \Omega^2(M)$, complex structure $\func{I}{\tang M}{\tang M}$, and associated Riemannian metric $g$.  Then $\mathcal{J}_\omega$ and $\mathcal{J}_I$ commute and
\[
  G:=-\mathcal{J}_{\omega} \circ \mathcal{J}_I
  = \begin{pmatrix} 0 & g^{\sharp} \\ g^{\flat} & 0 \end{pmatrix}
\]
is positive def\/inite, so $(M, \mathcal{J}_{\omega}, \mathcal{J}_I)$ is a GK manifold.
\end{enumerate}
\end{Example}

\begin{Example} \label{example2}
Let $(M_1,\mathcal{J}_1)$ and $(M_2,\mathcal{J}_2)$ be almost GC manifolds.  Then the direct sum~$\mathcal{J}$ of~$\mathcal{J}_1$ and~$\mathcal{J}_2$ is a map
  $
    \func{\mathcal{J}:=(\mathcal{J}_1, \mathcal{J}_2)}{
    \pont{M_1} \oplus \pont{M_2}
    }{
    \pont{M_1} \oplus \pont{M_2}
    }
  $,
which under the identif\/ication
  $
    \pont{M_1} \oplus \pont{M_2}
    \cong
    \pont{(M_1 \times M_2)}
  $
yields an almost GC structure on $M_1 \times M_2$.  We will call this the \deff{direct sum} of the almost GC structures on $M_1$ and $M_2$.  It is not hard to see that $(\mathcal{J}_1,\mathcal{J}_2)$ is a GC structure on $M_1 \times M_2$ if and only if $\mathcal{J}_i$ is a GC structure on $M_i$ for $i=1,2$.

Let $H_1 \in \Omega_{\text{cl}}^3(M_1)$ and $H_2 \in \Omega_{\text{cl}}^3(M_2)$, let $\func{\pi_i}{M_1 \times M_2}{M_i}$ be the natural projection for $i=1,2$, and set
$
  H := \pi_1^* H_1 + \pi_2^* H_2
$.
By the naturality of the exterior derivative, we know $H$ is a~closed three-form on $M_1 \times M_2$.  Furthermore, it is not hard to see that $(\mathcal{J}_1,\mathcal{J}_2)$ is an $H$-twisted GC structure on $M_1 \times M_2$ if and only if $\mathcal{J}_i$ is an $H_i$-twisted GC structure on $M_i$ for $i=1,2$.

There is a completely analogous product construction for almost GK and GK  manifolds as well.
\end{Example}

Let $(M,E,H)$ be a twisted GC manifold.  Suppose $S$ is a submanifold of $M$ given by the embedding $\hookfunc{j}{S}{M}$. Although $j$ induces a natural embedding $\hookfunc{j_*}{\tang S}{\tang M}$ of tangent bundles, because of the contravariance of cotangent bundles there is in general no obvious embedding $\pont{S} \hookrightarrow \pont{M}$ of the Pontryagin bundles.  The following def\/inition comes from~\cite{benbassat}.

For each $x \in S$, def\/ine
  \[
    E_{S,x}
    :=
    \left\{ \left( X, \lambda|_{S} \right) \in
    \mathbb{T}_{\mathbb{C},x} S
    \mid (X, \lambda) \in (\tang_{\mathbb{C},x} S \oplus \tang_{\mathbb{C},x}^* M ) \cap E_x \right\}.
  \]
By \cite[Lemma 8.2]{benbassat}, this $E_{S,x}$ is a maximally isotropic complex subspace of $\mathbb{T}_{\mathbb{C}, x} S$.
Let $E_S := \bigsqcup_{x \in S} E_{S, x}$.  Then $E_S$ is a constant-rank complex linear distribution of $\cpont{S}$, but is not in general a smooth subbundle, nor will it generally satisfy $E_S \cap \overline{E_S} = 0$.

\begin{Proposition} \label{prop6}
Let $(M,E,H)$ be a twisted GC manifold, let $\hookfunc{j}{S}{M}$ be a submanifold, and let $E_S \subset \cpont{S}$ be as defined above.  If $E_S$ is a subbundle of $\cpont{S}$, then $E_S$ is ($j^*H$)-twisted Courant involutive.
\end{Proposition}

In the untwisted case, where $H = 0$, Proposition~\ref{prop6} was proved in \cite[Corollary 3.1.4]{courant}.  The proof in the twisted case is nearly identical, with only minor changes to this proof and the relevant def\/initions and precursory results, (i.e.\ Def\/inition~2.3.2, Propositions~2.3.3 and~3.1.3, and Corollary~3.1.4 in~\cite{courant}).

\begin{Definition}
Let $(M,E,H)$ be a twisted GC manifold, and let $\hookfunc{j}{S}{M}$ be a submanifold.  If $E_S \subset \cpont{S}$ is a subbundle and satisf\/ies $E_S \cap \overline{E_S} = 0$, then $(S,E_S,j^* H)$ is a \deff{(twisted) GC submanifold} of $(M,E,H)$, and we denote by $\mathcal{J}_S$ the GC structure on $S$ induced by $E_S$.
\end{Definition}

\begin{Remark} \label{rem2}
Suppose $(M,E,H)$ is a twisted GC manifold, and $\hookfunc{j}{S}{M}$ is an \emph{open} submanifold.  Then since $\tang S$ and $\tang^* S$ can be identif\/ied with $(\tang M)|_S$ and $(\tang^* M)|_S$, respectively, we see that we can identify $E_S$ with $E|_S$, that $\mathcal{J}_S = \mathcal{J}|_{\pont{S}}$, and that $j^* H$ can be identif\/ied with $H|_S$.  Therefore an open submanifold of an $H$-twisted GC manifold is automatically an $H$-twisted GC manifold.  Similarly, an open submanifold of an $H$-twisted GK manifold is automatically an $H$-twisted GK manifold.
\end{Remark}

\begin{Definition}
Let $(M, \mathcal{J})$ be an almost GC manifold, and let $S \subset M$ be a submanifold.  A \deff{splitting bundle} for $S$ with respect to $(M, \mathcal{J})$ is a subbundle $N$ of $\tang M|_S \to S$ such that $\tang M|_S = \tang S \oplus N$ and $\tang S \oplus \operatorname{Ann} (N) \subset \pont{M}$ is invariant under $\mathcal{J}$. If a splitting bundle exists for~$S$, then $S$ is called a \deff{split submanifold} of $(M,\mathcal{J})$.
\end{Definition}

The following is an extension of Proposition 5.12 of \cite{benbassat} to the twisted case.  As with Proposition~\ref{prop6}, the original proof still holds with only minor alterations.

\begin{Proposition}
Let $(M,\mathcal{J},H)$ be a twisted GC manifold, and
let $\hookfunc{i}{S}{M}$ be a split submanifold of $M$ with splitting bundle $N \to S$.  Then $S$ is an $(i^* H)$-twisted GC submanifold of $(M,\mathcal{J},H)$, and the GC structure corresponding to the bundle $E_S$ is the same as the one induced by the restriction of $\mathcal{J}$ via the natural isomorphism
\[
  \pont{S} \cong \tang S \oplus \operatorname{Ann} N \subset \pont{M}.
\]
\end{Proposition}

It is straightforward to show that this implies the following.

\begin{Corollary}
Let $(M, \mathcal{J}_1, \mathcal{J}_2,H)$ be a twisted GK manifold, and and
let $\hookfunc{i}{S}{M}$ be a split submanifold of $M$ with respect to both $\mathcal{J}_1$ and $\mathcal{J}_2$, with common splitting bundle $N \to S$.  Then $\left( S, (\mathcal{J}_1)_W, (\mathcal{J}_2)_W \right)$ is an $(i^*H)$-twisted GK manifold.
\end{Corollary}

\begin{Definition}
Let $M$ be a manifold, and let $G$ be a Lie group acting smoothly on $M$.  This lifts to an action of $G$ on $\pont{M}$ by bundle automorphisms, given by
\[
  \func{\left( g_*, (g^{-1})^* \right)}{\pont{M}}{\pont{M}}
\]
for each $g \in G$, where $g_{\ast}$ is the {\bf pushforward} of tangent vectors by the map $\func{g}{M}{M}$ and $(g^{-1})^{\ast}$ is the {\bf pullback} of tangent covectors by the map $\func{g^{-1}}{M}{M}$.

Let $\mathcal{J}$ be an $H$-twisted GC structure on $M$.  We say that the $G$-action on $(M, \mathcal{J}, H)$ is \deff{canonical} if the following hold.
\begin{enumerate}\itemsep=0pt
  \item The dif\/ferential form $H$ is $G$-invariant, i.e.\ $g^* H = H$ for all $g \in G$.
  \item The action of $G$ on $\pont{M}$ commutes with $\mathcal{J}$, i.e.\ the diagram
\[
  \xymatrix{
  \pont{M} \ar[d]_{ \left( g_*, (g^{-1})^* \right) } \ar[r]^{\mathcal{J}}
  &
  \pont{M} \ar[d]^{ \left( g_*, (g^{-1})^* \right) } \\
  \pont{M} \ar[r]_{\mathcal{J}}
  &
  \pont{M}
  }
\]
commutes for all $g \in G$.
\end{enumerate}
\end{Definition}

It is easy to check that a smooth group action on a manifold commutes with an almost GC structure $\func{\mathcal{J}}{\pont{M}}{\pont{M}}$ if and only if the complex linear extension of the action preserves the corresponding complex Dirac structure.

\begin{Example}
Let $(M,\omega)$ be an almost symplectic manifold, let $\func{\omega^{\flat}}{\tang M}{\tang^* M}$ be the associated bundle isomorphism, and let
\[
  \mathcal{J}_{\omega} :=
  \begin{pmatrix}
    0 & - \omega^{\sharp} \\
    \omega^{\flat} & 0
  \end{pmatrix}
\]
be the associated almost GC structure on $M$.  Let $G$ be a Lie group acting smoothly on $M$.  It is easy to check that the $G$-action is symplectic if and only if the map $\func{\omega^{\flat}}{\tang M}{\tang^* M}$ is $G$-equivariant with respect to the pushforward action on $\tang M$ and the inverse pullback action on $\tang^* M$, if and only if the $G$-action commutes with $\mathcal{J}_{\omega}$.
\end{Example}

Recall that for a smooth action of a compact Lie group $G$ on a manifold~$M$, each connected component of the f\/ixed point set $M^G$ is a closed submanifold of $M$.  (Dif\/ferent components of~$M^G$ may have dif\/ferent dimensions.)

\begin{Proposition}  \label{prop9}
Let $M$ be a manifold, and let $\mathcal{J}$ be an almost GC structure on $M$. Suppose the compact Lie group $G$ acts canonically on $(M,\mathcal{J})$.  Then each component of $M^{G}$ is a split submanifold of $(M,\mathcal{J})$.
\end{Proposition}

\begin{proof}
Let $(M^G)'$ be a component of $M^G$.
First, recall that for each $x \in (M^G)'$ the derivative of the action of $G$ at $x$ def\/ines a linear action of $G$ on $\tang_x M$, and that $\tang_x (M^G)' = (\tang_x M)^G$.  Let $dg$ be a bi-invariant Haar measure on $G$, adjusted so that $dg(G) = 1$.  Def\/ine a bundle map $\func{\pi}{(\tang M)|_{(M^G)'}}{\tang (M^G)'}$ by setting
\[
  \pi_x(v) := \int_G (g \cdot v) \, dg
  \qquad
  \text{for all} \quad v \in V,
\]
for each $x \in (M^G)'$.  Def\/ine the subbundle $N \subset (\tang M)|_{(M^G)'}$ by setting $N_x := \ker \pi_x$ for each $x \in \left( M^G \right)'$.  That $\pi$ is a bundle map and $N$ is a vector bundle follow from the naturality of the technique of averaging by integration.  Note also that $\tang_x M = (\tang_x M)^G \oplus N_x = \tang_x (M^G)' \oplus N_x$ for each $x \in \left( M^G \right)'$.  To conclude that $N$ is a splitting for $(M^G)' \subset (M,\mathcal{J})$, it remains only to show that $\tang (M^G)' \oplus \operatorname{Ann}(N) = (\tang M)^G|_{(M^G)'} \oplus \operatorname{Ann}(N)$ is preserved by $\mathcal{J}$.

Fix $x \in \left( M^G \right)'$, let $V = \tang_x M$, and let $W = N_x$.  We claim that $\operatorname{Ann}(W) = (V^*)^G$, which would imply that $(V \oplus V^*)^G = V^G \oplus (V^*)^G = V^G \oplus \operatorname{Ann}(W)$.  Since $\mathcal{J}_x$ commutes with the action of $G$ on $V \oplus V^*$, we know that $\mathcal{J} \left( (V \oplus V^*)^G \right) = (V \oplus V^*)^G$, and hence a proof of this claim completes the proof of this proposition.

Let $\lambda \in \operatorname{Ann}(W)$, let $g \in G$, and let $u \in V$.  Decompose $u$ as $u=v + w$ for $v \in V^G$ and $w \in W$.
Then
\begin{alignat*}{2}
  (g \cdot \lambda)(u)
  &=
  \lambda ( g^{-1} \cdot u)
  && \\
  &=
  \lambda ( g^{-1} \cdot v) + \lambda(g^{-1} \cdot w)
  && \\
  &=
  \lambda (v) + \lambda(g^{-1} \cdot w)
  && \text{since } v \in V^G \\
  &=
  \lambda (v) + 0
  && \text{since $w \in W$, $W$ is $G$-stable, and }
  \lambda \in \operatorname{Ann}(W) \\
  &=
  \lambda (v) + \lambda (w)
  && \text{since $w \in W$ and } \lambda \in \operatorname{Ann}(W) \\
  &=
  \lambda (u).
  &&
\end{alignat*}
Therefore $\lambda \in (V^*)^G$, so $\operatorname{Ann}(W) \subset (V^*)^G$.  Note that
\[
  \operatorname{Ann}(W)
  \cong
  (V / W)^*
  =
  \left( (V^G \oplus W) / W \right)^*
  \cong
  (V^G)^*.
\]
Hence $\dim \operatorname{Ann}(W) = \dim (V^G)$, so $\operatorname{Ann} (W) = (V^*)^G$.
\end{proof}

\section{Background information on $G$-spaces}\label{section2}

In this section we give some brief def\/initions and results about compact group actions on manifolds which will be required in later sections.  The standard reference for the material on equivariant cohomology is \cite{supersymmetry}.  The material on orbit spaces and their stratif\/ication by orbit types can be found in \cite[Chapter 2]{duis} and \cite[Chapter 2]{ortega-ratiu}.

\subsection{Equivariant cohomology}

Let $M$ be a manifold and $G$ be a compact Lie group acting smoothly on $M$.  Consider the space $\Omega^{k}(M) \otimes S^i(\mathfrak{g}^*)$, where $S^i$ denotes the degree $i$ elements of the symmetric algebra.  This is a~$G$-space with action def\/ined by linear extension of the rule $g \cdot (\alpha \otimes p) := \left( (g^{-1})^* \alpha \right) \otimes (p \circ \mathrm{Ad}_{g^{-1}})$ for $g \in G$, $\alpha \in \Omega^{\star}(M)$, $p \in S(\mathfrak{g}^*)$.  We can identify $\Omega^{k}(M) \otimes S^i(\mathfrak{g}^*)$ with the space of degree $i$ polynomial maps $\mathfrak{g} \to \Omega^k(M)$ via
\[
  \alpha \otimes p \colon \xi \mapsto p(\xi) \cdot \alpha
\]
for $\xi \in \mathfrak{g}$.  An element of $\Omega^{k}(M) \otimes S^i(\mathfrak{g}^*)$ is $G$-invariant if and only if its corresponding polynomial map is $G$-equivariant with respect to the adjoint action of $G$ on $\mathfrak{g}$ and the action of $G$ on $\Omega^k(M)$ given by $g \cdot \alpha := (g^{-1})^* \alpha$ for $g \in G$, $\alpha \in \Omega^k(M)$.

\begin{Definition}
Let $M$ be a manifold and $G$ be a compact Lie group acting smoothly on $M$.  The space of \deff{equivariant differential forms of degree $n$} on $M$ is
\[
  \Omega_G^n(M) := \bigoplus_{i=0}^{\lfloor n/2 \rfloor}
  \left( \Omega^{n-2i}(M) \otimes S^i(\mathfrak{g}^*) \right)^G.
\]
The dif\/ferential $\func{\mathrm{d}_G}{\Omega_G^n}{\Omega_G^{n+1}}$ is def\/ined, viewing equivariant forms as maps $\mathfrak{g} \to \Omega^{\star}(M)$, by
\[
  \mathrm{d}_G ( \alpha \otimes p )(\xi)
  :=
  \left( \mathrm{d} \alpha - \iota_{\xi_M} \alpha \right) \, p(\xi)
  \quad \text{for all } \xi \in \mathfrak{g}.
\]
The \deff{Cartan model for the $G$-equivariant cohomology} of $M$ is $H_G^{\star}(M) := H^{\star} (\Omega_G^{\star}, \mathrm{d}_G)$.
\end{Definition}

Suppose now that $G$ acts freely on $M$.  Then the $G$-equivariant cohomology of $M$ is naturally isomorphic as a graded algebra to the de Rham cohomology of the quotient $M/G$,
\[
  H_G^{\star}(M) \cong H^{\star} (M/G).
\]
We denote this isomorphism by $\func{\kappa}{H_G^{\star}(M)}{H^{\star}(M/G)}$.

Let $B \in \Omega^n(M)$.  The form $B$ is called \deff{basic} if it is $G$-equivariant and if $\iota_{\xi_M} B = 0$ for all $\xi \in \mathfrak{g}$.  If there is a dif\/ferential form $\widetilde{B} \in \Omega^n(M/G)$ such that the pullback of $\widetilde{B}$ by the quotient map $M \to M/G$ equals $B$, then we say $B$ \deff{descends} to $\widetilde{B}$.

\begin{Proposition}
Let $M$ be a manifold and $G$ be a compact Lie group acting smoothly and freely on $M$.
\begin{enumerate} \itemsep=0pt\litlet
  \item If $B \in \Omega^n(M)$ is basic, then $B$ descends to some $\widetilde{B} \in \Omega^n(M/G)$.
  \item If $B \in \Omega^n(M)^G \subset \Omega_G^n(M)$ is equivariantly closed, i.e.\ $\mathrm{d}_G B = 0$, then $B$ is closed and basic and descends to some closed $\widetilde{B} \in \Omega^n(M/G)$ such that
    \[
      \kappa[B] = [\widetilde{B}],
    \]
  where $[B]$ and $[\widetilde{B}]$ are the cohomology classes of $B$ and $\widetilde{B}$, respectively.
  \item If $\eta \in \Omega_G^n(M)$ is equivariantly closed, then there exists $\Gamma \in \Omega_G^{n-1}(M)$ so that $\eta + \mathrm{d}_G \Gamma \in \Omega^n(M)^G \subset \Omega_G^n(M)$.  In this case, since $\eta + \mathrm{d}_G \Gamma$ is equivariantly closed, it descends to some $\widetilde{\eta} \in \Omega^n(M/G)$ such that $\kappa[\eta] = [\widetilde{\eta}]$.
\end{enumerate}
\end{Proposition}

\begin{Definition} \label{def3}
Let $M$ be a manifold and $G$ be a compact Lie group acting on $M$ smoothly and freely.  Then $M \to M/G$ is a (left) principal $G$-bundle.  A \deff{connection} on this bundle is a~$\mathfrak{g}$-valued one-form $\theta \in \Omega^1(M,\mathfrak{g})$ such that
\begin{enumerate}\itemsep=0pt
  \item $\theta$ is $G$-equivariant, i.e.\ $g^* \theta = \mathrm{Ad}_g \circ \theta$;
  \item $\theta (\xi_M) \equiv \xi$ for all $\xi \in \mathfrak{g}$.
\end{enumerate}
\end{Definition}

\subsection{Orbit type stratif\/ication}

Let $G$ be a group.  For each subgroup $H$ of $G$, we will denote by $(H)$ the set of subgroups of~$G$ that are conjugate to $H$.  Suppose $G$ is a compact Lie group and $M$ is a manifold on which~$G$ acts smoothly.  Note that the conjugacy relation among subgroups of~$G$ preserves closedness, and hence also preserves the property of being a Lie subgroup.

\begin{Definition}
Let $x \in M$, and let $G_x := \{ g \in G \mid g \cdot x = x \}$ be the isotropy subgroup of~$x$ in~$G$.  The \deff{orbit type} of the point $x$, or of $G \cdot x$, is the set $(G_x)$ of subgroups of $G$ that are conjugate to $G_x$.

Let $H$ be a closed subgroup $H$ of $G$.  The \deff{$(H)$-orbit type submanifold} of $M$ is the set $M_{(H)} := \{ x \in M \mid G_x \in (H) \}$.  The \deff{$H$-isotropy type submanifold} of $M$ is the set $M_H := \{ x \in M \mid G_x = H \}$.  The \deff{$H$-fixed point submanifold} of $M$ is the set $M^H := \{ x \in M \mid G_x \subset H \}$.
\end{Definition}

Note that the sets def\/ined above are related by the equation $M_H = M_{(H)} \cap M^H$.
Also, two $G$-orbits in $M$ have the same orbit type if and only if they are $G$-equivariantly dif\/feomorphic.  This leads one to the following def\/initions.

\begin{Definition}
Let $x \in M$, and let $H = G_x$.  The \deff{local action type submanifold through~$x$} is the subset $M_{(H)}^{l_x} \subset M$ of points $y \in M$ such that there is a $G$-equivariant dif\/feomorphism between $G$-invariant open neighborhoods of $x$ and $y$.  Def\/ine $M_H^{l_x} := M_{(H)}^{l_x} \cap M^H$.
\end{Definition}

Some important properties of the sets we have def\/ined above are collected in the following proposition.  Their proofs can be found in the references cited at the beginning of this section.

\begin{Proposition} \label{prop2}
Let $G$ be a compact Lie group, and $M$ be a manifold on which $G$ acts smoothly.  Let $x \in M$ and put $H = G_x$.  Then the following hold.
\begin{enumerate} \itemsep=0pt\litlet

  \item $M_{(H)}^{l_x}$, respectively $M_{H}^{l_x}$, is an open and closed subset of $M_{(H)}$, respectively $M_H$.

  \item The sets $M_{(H)}^{l_x}$ and $M_H^{l_x}$ are locally closed embedded submanifolds of $M$, as is each connected component of $M_H$, of $M_{(H)}$, and of $M^H$.

  \item $M_H^{l_x}$ and $M_{(H)}^{l_x}$ consists of the union of certain components of $M_H$ and $M_{(H)}$, respectively.

  \item $M_H$ and $M_H^{l_x}$ are open in $M^H$.

  \item $M_{(H)}$ and $M_{(H)}^{l_x}$ are $G$-stable, and $G \cdot M_H = M_{(H)}$ and $G \cdot M_H^{l_x} = M_{(H)}^{l_x}$.

  \item Let $N=N_G(H)$ be the normalizer of $H$ in $G$. Both $M_H$ and $M_H^{l_x}$ are $N$-stable, and $N/H$ acts freely on both.  Hence $M_H^{l_x} / N \cong M_H^{l_x} / \left( N/H \right)$ is a manifold.

  \item The inclusions $M_H \hookrightarrow M_{(H)}$ and $M_H^{l_x} \hookrightarrow M_{(H)}^{l_x}$ induce homeomorphisms $M_H \to M_{(H)}$ and $M_H^{l_x} / N \to M_{(H)}^{l_x} / G$.  Thus the quotient $M_{(H)}^{l_x} / G$ inherits a natural manifold structure.

  \item Each component of the quotient $M_{(H)} / G$ inherits a natural manifold structure.
\end{enumerate}
\end{Proposition}

In general, the orbit space $M/G$ can be a very singular space.  It will be a Hausdorf\/f and second-countable topological space, but will rarely inherit a manifold, or even an orbifold, structure from $M$.
However, because $M$ is the disjoint union of its orbit type submanifolds, we can also partition the orbit space:
\begin{equation}
 \label{eq5}
   M/G = \bigsqcup_{(H)} M_{(H)}/G,
\end{equation}
where the disjoint union is taken over all the distinct orbit type submanifolds of $M$.
Since each component of $M_{(H)}/G$ is a manifold, we know that, after ref\/ining the partition to components, (\ref{eq5}) is a partition of $M/G$ into manifolds.  It is called the \deff{orbit type partition} of $M/G$.

\begin{Remark}
All of the above results hold true even if $G$ is an \emph{arbitrary} Lie group, so long as it acts on $M$ both smoothly and \emph{properly}.
\end{Remark}

\section{Hamiltonian actions on generalized complex manifolds}\label{section3}

In \cite{lintolman}, the authors proposed the following def\/inition of Hamiltonian actions on GC manifolds.

\begin{Definition}
\label{genhamactiondef}

Let $(M,\mathcal{J})$ be an untwisted GC manifold, let $E$ be the associated complex Dirac structure on $M$, and let $G$ be a Lie group acting canonically on $(M,\mathcal{J})$.  This action is \deff{generalized Hamiltonian} if there exists a $G$-equivariant map $\func{\mu}{M}{\mathfrak{g}^{\ast}}$ such that, for all $\xi \in \mathfrak{g}$,
\[
  \xi_M = - \mathcal{J} ( \mathrm{d} \mu^{\xi} )
\]
or equivalently $\xi_M - i \mathrm{d} \mu^{\xi} \in \Gamma(E)$.
Here $\func{\mu^{\xi}}{M}{\mathbb{R}}$ is the smooth function def\/ined by $\mu^{\xi}(x) := \iprod{\mu(x)}{\xi}$ for all $x \in M$.  The map $\mu$ is called a \deff{generalized moment map} for the $G$-action on~$(M, \mathcal{J})$.

Let $(M,\mathcal{J},H)$ be a twisted GC manifold, and let $G$ be a Lie group acting canonically on $(M,\mathcal{J},H)$.  This action is \deff{twisted generalized Hamiltonian} if there exists a $G$-equivariant map $\func{\mu}{M}{\mathfrak{g}^{\ast}}$ and a $G$-equivariant $\mathfrak{g}^*$-valued one-form $\alpha \in \Omega^1(M, \mathfrak{g}^{\ast})$ on $M$ such that, for all $\xi \in M$,
\begin{enumerate}\itemsep=0pt
  \item $\xi_M = - \mathcal{J} ( \mathrm{d} \mu^{\xi} ) - \alpha^{\xi}$ (or equivalently $\xi_M + \alpha^{\xi} - i \mathrm{d} \mu^{\xi} \in \Gamma(E)$), and
  \item $\iota_{\xi_M} H = \mathrm{d} \alpha^{\xi}$.
\end{enumerate}
Here $\mu^{\xi}$ is as def\/ined above, and $\alpha^{\xi} \in \Omega^1(M)$ is the dif\/ferential one-form on $M$ def\/ined by $(\alpha^{\xi})_x (v) := \iprod{\alpha_x(v)}{\xi}$ for all $x \in M$, $v \in \tang_x M$.  The map $\mu$ and the one-form $\alpha$ are called a~\deff{generalized moment map} and a \deff{moment one-form}, respectively, for the action $G$-action on $(M, \mathcal{J}, H)$.
\end{Definition}

\begin{Definition}
Let $(M,\mathcal{J}_1,\mathcal{J}_2)$ be a GK manifold, and let $G$ be a Lie group acting on $M$ and preserving both $\mathcal{J}_1$ and $\mathcal{J}_2$.  This action is called \deff{generalized Hamiltonian} if the action of $G$ on $(M,\mathcal{J}_1)$ is generalized Hamiltonian.

Similarly, if $(M,\mathcal{J}_1,\mathcal{J}_2,H)$ be a twisted GK manifold, and the $G$-action preserves $\mathcal{J}_1$, $\mathcal{J}_2$, and $H$, then the action is \deff{twisted generalized Hamiltonian} if the action of $G$ on $(M,\mathcal{J}_1,H)$ is twisted generalized Hamiltonian.
\end{Definition}

\begin{Remark} \label{rem1} \qquad
\begin{enumerate} \itemsep=0pt\litlet
  \item Note that a moment one-form $\alpha \in \Omega^1(M,\mathfrak{g}^*)$ is an equivariant dif\/ferential form of degree~$3$.

  \item Because $E$ is an isotropic subbundle, the condition that $\xi_M + \alpha^{\xi} - i \mathrm{d} \mu^{\xi} \in E$ implies that $\iiprod{\xi_M + \alpha^{\xi} - i \mathrm{d} \mu^{\xi}}{\xi_M + \alpha^{\xi} - i \mathrm{d} \mu^{\xi}}=0$, and hence that $\iota_{\xi_M} \alpha^{\xi} = \iota_{\xi_M} \mathrm{d} \mu^{\xi} = 0$.

  \item Given a GC manifold $(M,\mathcal{J})$, one can consider this as a twisted GC manifold $(M,\mathcal{J},H)$ by setting $H=0$.  Therefore, if a Lie group $G$ acts on $(M,\mathcal{J})$ canonically, we have two notions of whether the action is Hamiltonian.  It may be Hamiltonian as an action on~$(M,\mathcal{J})$, in which case there is just a moment map, or it may be Hamiltonian as an action on $(M,\mathcal{J},H)$, in which case there is {\bf both} a moment map {\bf and} a moment one-form.
  It is potentially interesting to explore both possibilities.
\end{enumerate}
\end{Remark}

\begin{Example} \label{example3}
Let $(M,\omega)$ be a symplectic manifold, and let $G$ be a Lie group acting on $(M,\omega)$ in a Hamiltonian fashion with moment map $\func{\Phi}{M}{\mathfrak{g}^*}$.  Recall that this means the $G$-action is symplectic, the map $\Phi$ is $G$-equivariant, and for all $\xi \in \mathfrak{g}$ we have $\mathrm{d} \Phi^{\xi} = \iota_{\xi_M}$.
Let $\mathcal{J}_{\omega}$ be the GC structure on $M$ induced by $\omega$.  As discussed in Example 3.8 of \cite{lintolman}, the action of $G$ on $(M,\mathcal{J}_{\omega})$ is generalized Hamiltonian, and $\Phi$ is a generalized moment map.
\end{Example}

\begin{Theorem} \label{thm6}
Let $(M,E,H)$ be a twisted GC manifold, where $E$ is the associated complex Dirac structure, and let $G$ be a Lie group acting on $(M,E,H)$ in a Hamiltonian fashion with moment map $\func{\mu}{M}{\mathfrak{g}^*}$ and moment one-form $\alpha \in \Omega^1(M,\mathfrak{g}^*)$.
If $\hookfunc{j}{S}{M}$ is a $G$-stable twisted GC submanifold of $(M,E,H)$, then the restriction of the action of $G$ to $(S,E_S,j^* H)$ is Hamiltonian with moment map $\func{\mu|_S}{S}{\mathfrak{g}^*}$ and moment one-form $j^* \alpha \in \Omega^1(S,\mathfrak{g}^*)$.
\end{Theorem}

\begin{proof}
First we will prove that the action of $G$ on $S$ preserves $E_S$.  Let $x \in S$ and $(X, \lambda) \in (\tang_{\mathbb{C}, x} S \oplus \tang_{\mathbb{C}, x}^{\ast} M) \cap E_x$, which means that $(X, \lambda|_S) \in E_{S,x}$.  Then for any $g \in G$ we have
\[
  g \cdot (X + j^* \lambda) = g_* (X) + (g^{-1})^* \lambda|_S.
\]
Because $S$ is $G$-stable, the inclusion $\hookfunc{j}{S}{M}$ is $G$-equivariant, i.e. the $G$-action commutes with $j$.  Hence $\hookfunc{j_*}{\tang S}{\tang M}$ is $G$-equivariant, so $g_* (X) \in \tang_{g \cdot x} S$.  Also
\[
  (g^{-1})^* \left( \lambda|_S \right)
  =
  (g^{-1})^* j^* \lambda
  =
  j^* (g^{-1})^* \lambda
  =
  \left( (g^{-1})^* \lambda \right)|_S.
\]
Since $E$ is $G$-stable, we have $g \cdot (X + \lambda) = g_* (X) + (g^{-1})^* \lambda \in E_{g \cdot x}$.  Therefore $g \cdot (X, j^* \lambda) = \left( g_* (X), j^* (g^{-1})^* \lambda \right) \in E_{S, g \cdot x}$.  Thus $E_S$ is $G$-stable.

Now suppose that $(S,E_S,j^* H)$ is a GC submanifold of $(M,E,H)$, meaning that $E_S$ is a~vector bundle, that $E_S \cap \overline{E_S} = 0$, and that $E_S$ is $j^* H$-twisted Courant involutive.  Since $j$ is $G$-equivariant, for all $\xi \in \mathfrak{g}$ we have $\xi_M|_S = \xi_S$, $(j^* \mu)^{\xi} = j^* (\mu^{\xi})$, and $(j^* \alpha)^{\xi} = j^* (\alpha^{\xi})$.  Furthermore, by the naturality of the exterior derivative we have
\[
  \mathrm{d} (j^* \mu)^{\xi}
  =
  \mathrm{d} j^* (\mu^{\xi})
  =
  j^* \big(\mathrm{d} \mu^{\xi}\big),
\]
so $\mathrm{d} (\mu|_{S})^{\xi} = (\mathrm{d} \mu^{\xi})|_S$.  For each $x \in S \subset M$, since $\eval{\left( \xi_M + \alpha^{\xi} - i \, \mathrm{d} \mu^{\xi} \right)}{x} \in E_x$, this means that
\[
  \eval{\big( \xi_S + \alpha^{\xi}|_S - i \, \mathrm{d} \mu^{\xi}\big|_S \big)}{x} \in E_{S,x}.
\]
Again using the $G$-equivariance of $j$, for all $x \in S$ we have
\[
  \eval{\iota_{\xi_S} j^* (H)}{x}
  =
  \eval{j^* (\iota_{\xi_S} H)}{x}
  =
  \eval{j^* (\iota_{\xi_M} H)}{x}
  =
  \eval{j^* (\mathrm{d} \alpha^{\xi})}{x}
  =
  \eval{\mathrm{d} (j^* \alpha)^{\xi}}{x}.
\]
Thus the action of $G$ on $(S,E_S,j^* H)$ is twisted Hamiltonian with moment map $\mu|_S$ and moment one-form $\alpha|_S$.
\end{proof}

The above result holds also for the untwisted case, of course, by putting $H=0$ and $\alpha=0$.

The following three results are exactly what makes reduction of generalized Hamiltonian manifolds possible.

\begin{Theorem}[Lemma 3.8 and Proposition 4.6 of \cite{lintolman}] \label{thm4}
Let a compact Lie group $G$ act on a~GC manifold $(M,\mathcal{J})$, respectfully a GK manifold $(M,\mathcal{J}_1,\mathcal{J}_2)$, in a Hamiltonian fashion with moment map $\func{\mu}{M}{\mathfrak{g}^*}$.  Suppose $a \in \mathfrak{g}^*$ is an element such that $G$ acts freely on the inverse image $\mu^{-1}(\mathcal{O}_a)$ of the coadjoint orbit $\mathcal{O}_a$ of $G$ through $a$.  Then the quotient space $\mu^{-1}(\mathcal{O}_a) / G$ inherits a natural GC structure $\widetilde{\mathcal{J}}$ from $\mathcal{J}$, respectfully a natural GK structure $(\widetilde{\mathcal{J}}_1,\widetilde{\mathcal{J}}_2)$ from~$(\mathcal{J}_1,\mathcal{J}_2)$.
\end{Theorem}

\begin{Lemma}[Lemma A.6 of \cite{lintolman}]
\label{lem2}
Let a compact Lie group $G$ act freely on a manifold $M$.  Let $H$ be a $G$-invariant and closed three-form, and let $\func{\alpha}{\mathfrak{g}}{\Omega^1(M)}$ be an equivariant map.  Fix a connection $\theta \in \Omega^1(M, \mathfrak{g})$ on the principal $G$-bundle $M \to M/G$.  Then if $H+\alpha \in \Omega_G^3(M)$ is equivariantly closed, there exists a natural form $\Gamma \in \Omega^2(M)^G$ so that $\iota_{\xi_M} \Gamma = \alpha^{\xi}$ for all $\xi \in \mathfrak{g}$.  Thus $H + \alpha + \mathrm{d}_G \Gamma \in \Omega^3(M)^G \subset \Omega_G^3(M)$ is closed and basic and so descends to a closed form $\widetilde{H} \in \Omega^3(M/G)$ so that $[\widetilde{H}] = \kappa[H+\alpha]$.
\end{Lemma}

\begin{Theorem}[Propositions A.7 and A.10 of \cite{lintolman}] \label{thm5}
Let a compact Lie group $G$ act on a~twisted GC manifold $(M,\mathcal{J},H)$, respectfully a twisted GK manifold $(M,\mathcal{J}_1,\mathcal{J}_2,H)$, in a Hamiltonian fashion with moment map $\func{\mu}{M}{\mathfrak{g}^*}$ and moment one-form $\alpha \in \Omega^1(M,\mathfrak{g})$.  Suppose $a \in \mathfrak{g}^*$ is an element such that $G$ acts freely on the inverse image $\mu^{-1}(\mathcal{O}_a)$ of the coadjoint orbit~$\mathcal{O}_a$ of~$G$ through~$a$.  Assume that $H+\alpha$ is equivariantly closed.
Given a connection on the principal $G$-bundle $\mu^{-1}(\mathcal{O}_a) \to \mu^{-1}(\mathcal{O}_a)/G$, the quotient space $\mu^{-1}(\mathcal{O}_a) / G$ inherits an $\widetilde{H}$-twisted GC structure $\widetilde{\mathcal{J}}$ from $\mathcal{J}$, respectfully an $\widetilde{H}$-twisted GK structure $(\widetilde{\mathcal{J}}_1,\widetilde{\mathcal{J}}_2)$ from $(\mathcal{J}_1,\mathcal{J}_2)$, where $\widetilde{H}$ is defined as in Lemma~{\rm \ref{lem2}} above.  Up to $B$-transform, these inherited structures are independent of our choice of connection.
\end{Theorem}

\begin{Definition}
The quotient space $\mu^{-1}(\mathcal{O}_a) / G$ in Theorems~\ref{thm4} and \ref{thm5} is called the \deff{gene\-ralized complex quotient} (or \deff{generalized K\"ahler quotient}, as applicable), or the \deff{Lin--Tolman quotient}, of $M$ by $G$ at level $a$.  We use the notation
\[
  M_a := \mu^{-1}(\mathcal{O}_a) / G.
\]
\end{Definition}

\begin{Remark}
As noted in Example 3.9 of \cite{lintolman}, in the context of the hypotheses of Theorem~\ref{thm4}, if the GC structure and moment map come from a \emph{symplectic} structure and moment map, then the GC structure on the quotient is exactly the one induced by the Marsden--Weinstein ssymplectic structure on the quotient.
\end{Remark}

The following result will be useful to us later.  Its proof follows trivially from the def\/initions of generalized and twisted generalized Hamiltonian actions.

\begin{Lemma} \label{lem6}
Let $(M,\mathcal{J},H)$ be a twisted GC manifold with a Hamiltonian action of a Lie group $G$, moment map $\func{\mu}{M}{\mathfrak{g}^*}$, and moment one-form $\alpha \in \Omega^1(M,\mathfrak{g}^*)$.
Let $K \subset G$ be a~Lie subgroup.
Then the induced action of $K$ on $(M,\mathcal{J},H)$ is also Hamiltonian, with generalized moment map and moment one-form the compositions of $\mu$ and $\alpha$, respectively, with the projection $\mathfrak{g}^* \twoheadrightarrow \mathfrak{k}^*$ dual to the inclusion $\mathfrak{k} \hookrightarrow \mathfrak{g}$:
  \[
    \xymatrix{
      M \ar[r]^{\mu} & \mathfrak{g}^* \ar@{>>}[r] & \mathfrak{k}^*
    },
    \qquad
    \xymatrix{
      \tang M \ar[r]^{\alpha}
      & \mathfrak{g}^* \ar@{>>}[r] & \mathfrak{k}^*
    }.
  \]
\end{Lemma}

\begin{Example} \label{example5}
Let $G$ be a Lie group, and let $(M_i,\mathcal{J}_i,H_i)$ be a twisted GC manifold on which $G$ acts in a Hamiltonian fashion with moment map $\func{\mu_i}{M_i}{\mathfrak{g}^*}$ and moment one-form $\alpha_i \in \Omega^1(M_i,\mathfrak{g}^*)$, for $i=1,2$.  Let $(M_1 \times M_2, \mathcal{J}, H)$ be the product of these two GC manifolds, as def\/ined in Example~\ref{example2}.  Recall that $\mathcal{J} = (\mathcal{J}_1,\mathcal{J}_2)$ and $H = \pi_1^* H_1 + \pi_2^* H_2$, where $\func{\pi_i}{M_1 \times M_2}{M_i}$ is the natural projection for $i=1,2$.  Def\/ine $\func{\mu}{M_1 \times M_2}{\mathfrak{g}^* \oplus \mathfrak{g}^*}$ and $\alpha \in \Omega^1(M_1 \times M_2, \mathfrak{g}^* \oplus \mathfrak{g}^*)$ by $\mu = \pi_1^* \mu_1 + \pi_2^* \mu_2$ and $\alpha = \pi_1^* \alpha_1 + \pi_2^* \alpha_2$.  It is easy to check that the action of $G \times G$ on $M_1 \times M_2$ is twisted generalized Hamiltonian with moment map $\mu$ and moment one-form $\alpha$.

Embedding $G$ diagonally in $G \times G$, we obtain a Hamiltonian action of $G$ on $M_1 \times M_2$.  The projection $\mathfrak{g}^* \oplus \mathfrak{g}^* \twoheadrightarrow \mathfrak{g}^*$ induced by this embedding is given by addition: $(\lambda_1, \lambda_2) \mapsto \lambda_1 + \lambda_2$, so a moment map and moment one-form for the $G$-action on $M_1 \times M_2$ is given by
\[
  M_1 \times M_2 \to \mathfrak{g}^*,
  \qquad
  (x_1,x_2) \mapsto \mu_1(x_1) + \mu_2(x_2)
\]
and
\[
  \tang M_1 \times \tang M_2
  \to \mathfrak{g}^*,
  \qquad
  (X_1,X_2) \mapsto \alpha_1(X_1) + \alpha_2(X_2),
\]
respectively.
\end{Example}

Perhaps the most important instance of the construction of Example~\ref{example5} is if we start with an arbitrary twisted generalized Hamiltonian $G$-manifold, $(M,\mathcal{J},H,\mu,\alpha)$, and let the second GC manifold be a coadjoint orbit $\mathcal{O}_{a}$ in $\mathfrak{g}^*$, where $a \in \mathfrak{g}^*$ is some f\/ixed element.  Let $\omega_{a}$ be the canonical symplectic structure on $\mathcal{O}_{a}$.  The action of $G$ on $\mathcal{O}_{a}$ is Hamiltonian in the symplectic sense, with moment map given by the inclusion $\mathcal{O}_{a} \hookrightarrow \mathfrak{g}^*$.  Using the symplectic structure $-\omega_{a}$ instead, the action is still Hamiltonian, but now the moment map is given by the negative inclusion $\mathcal{O}_{a} \to \mathfrak{g}^*$, $\lambda \mapsto -\lambda$.

As described in Examples \ref{example4} and \ref{example3}, the symplectic structure $-\omega_{a}$ induces a GC structure~$\mathcal{J}_{a}$ on~$\mathcal{O}_{a}$, and the $G$-action on $\mathcal{O}_{a}$ is generalized Hamiltonian with the same moment map.  Viewing $(\mathcal{O}_{a},\mathcal{J}_{a})$ as a twisted GC manifold where the twisting is by the zero three-form, the $G$-action is twisted generalized Hamiltonian with a constantly vanishing moment one-form.  Then the diagonal $G$-action on $M \times \mathcal{O}_{a}$ is twisted generalized Hamiltonian with moment map
\[
  \func{\mu'}{M \times \mathcal{O}_{a}}{\mathfrak{g}^*},
  \qquad (x,\lambda) \mapsto \mu(x) - \lambda
\]
and moment one-form
\[
  \func{\alpha'}{
  \tang M \times \tang (\mathcal{O}_{a})}{
  \mathfrak{g}^*},
  \qquad
  (X,Y) \mapsto \alpha(X).
\]
The reason this construction is important is that it is the basis of the \deff{shifting trick}.  If one wishes to reduce $M$ by $G$ at level $a \in \mathfrak{g}^*$, one can instead consider the reduction of $M \times \mathcal{O}_{a}$ by~$G$ at level~$0$, because
\[
  M_a \approx (M \times \mathcal{O}_{a})_0
\]
as topological spaces.
To see this, observe that $\mu^{-1}(\mathcal{O}_a)$ and $(\mu')^{-1}(0)$ are $G$-equivariantly homeomorphic via the maps
\[
  \mu^{-1}(\mathcal{O}_a) \to (\mu')^{-1}(0),
  \qquad
  x \mapsto \left(x, \mu(x) \right)
\]
and
\[
  (\mu')^{-1}(0) \to \mu^{-1}(\mathcal{O}_a),
  \qquad
  (x, \lambda) \mapsto x.
\]

\section{Partition of the generalized reduced space}\label{section4}

Let $M$ be a manifold, $G$ be a Lie group acting on $M$ smoothly, and $\func{\mu}{M}{\mathfrak{g}^*}$ a smooth, $G$-equivariant map.
Let $a \in \mathfrak{g}^*$.
By equivariance the pre-image $\mu^{-1}(\mathcal{O}_a)$ of the coadjoint orbit~$\mathcal{O}_a$ is preserved by $G$, and so we can consider the quotient space $\mu^{-1}(\mathcal{O}_a) / G$.
Let $M = \bigsqcup M_{(H)}$ be the orbit type partition of $M$.
Because each set $M_{(H)}$ is stable under $G$, each intersection $\mu^{-1}(\mathcal{O}_a) \cap M_{(H)}$ is also stable under $G$,
so the orbit type partition of $M$ descends to a partition
\[
  \mu^{-1}(\mathcal{O}_a) / G
  = \bigsqcup_{(H)} \left( \mu^{-1}(\mathcal{O}_a) \cap M_{(H)} \right) / G
\]
of the quotient $\mu^{-1}(\mathcal{O}_a) / G$.

Suppose now $M$ is a symplectic manifold, the $G$-action is Hamiltonian, and $\mu$ is a moment map.  In this case the quotient space $M_a := \mu^{-1}(\mathcal{O}_a) / G$ is called the \deff{symplectic reduction}, or \deff{Marsden--Weinstein quotient}, of $M$ at level $a$.
The symplectic moment map condition is that $\mathrm{d} \mu^{\xi} = \iota_{\xi_M} \omega$ for all $\xi \in \mathfrak{g}$.  If $G$ acts freely on $\mu^{-1}(\mathcal{O}_a)$, then each $\xi_M$ is nonzero on $\mu^{-1}(\mathcal{O}_a)$, which by the non-degeneracy of $\omega$ implies that $a$ is a regular value of $\mu$.  Therefore $\mu^{-1}(\mathcal{O}_a) \subset M$ is a submanifold, so $M_a$ is a manifold.  In this case, Marsden and Weinstein proved that $M_a$ inherits a natural symplectic structure.  Theorems~\ref{thm4} and \ref{thm5}, proved in \cite{lintolman}, are analogues of this result.

In the event that the symplectic quotient is singular, one can consider the individual parts of the partitioned quotient.
In \cite{stratified}, Lerman and Sjamaar proved that each component of $\left( M_{a} \right)_{(H)} := \left( \mu^{-1}(\mathcal{O}_a) \cap M_{(H)} \right) / G$ inherits a natural symplectic structure.
The main results of this paper are analogues of this in the generalized complex case.

\begin{Remark} \label{rem3}
By the symplectic moment map condition, $\mathrm{d} \mu^{\xi} = \iota_{\xi_M} \omega$, if $a \in \mathfrak{g}^*$ is a regular value of $\mu$, then each vector f\/ield $\xi_M$ is nowhere zero on $\mu^{-1}(a)$.  This means that the action of~$G$ on $\mu^{-1}(a)$ is at least \emph{locally free}, which means that the quotient $M_a$ is at worst an orbifold, to which Marsden and Weinstein were able to associate a symplectic structure.  By Sard's Theorem, a generic value of $\mu$ will be regular, so the generic result of symplectic reduction is a symplectic orbifold.

If $(M,\mathcal{J})$ is an untwisted GC manifold with moment map $\mu$, then the generalized moment map condition, $\mathcal{J}(\mathrm{d} \mu^{\xi}) = -\xi_M$, likewise guarantees the equivalence of regular values and local freeness of the action, so $M_a$ is at worst an orbifold.  However, if $(M,\mathcal{J},H)$ is a twisted GC manifold with moment map $\mu$ and moment one-form $\alpha$, then this equivalence may no longer hold, due to the presence of the moment one-form in the moment condition:
\[
  \mathcal{J}(\mathrm{d} \mu^{\xi}) = -\xi_M - \alpha^{\xi}.
\]
Specif\/ically, $\xi_M$ could vanish even if $\mathcal{J}(\mathrm{d} \mu^{\xi})$ does not.  Therefore, it seems that the generic result of GC reduction may be a GC singular space.
\end{Remark}

Before stating and proving our main theorem, we need the following lemma.

\begin{Lemma} \label{lem1}
Let $(M,\mathcal{J},H)$ be a \emph{compact}, twisted GC manifold, and let $G$ be a compact Lie group acting on $(M,\mathcal{J},H)$ in a Hamiltonian fashion with moment map $\func{\mu}{M}{\mathfrak{g}^*}$ and moment one-form $\alpha \in \Omega^1(M,\mathfrak{g})$.
If the $G$-action on $M$ is trivial, then $\mathrm{d} \mu = \alpha \equiv 0$.
\end{Lemma}

In Lemma 5.5 of \cite{linbaird}, the authors proved the above result in the case that $G$ is a torus;  however, their proof holds just as well in the non-abelian case.  It relies on viewing the components of $\mu$ as the real parts of a pseudo-holomorphic function and applying a version of the Maximum Principle, a course f\/irst taken in \cite{nitta}.  A thorough description of this version of the Maximum Principle can be found in Section 4.4 of \cite{MyThesis}.

\begin{Theorem}[Singular generalized reduction] \label{MainThm} \qquad
\begin{enumerate} \itemsep=0pt\litlet
  \item Let $(M,\mathcal{J})$ be a GC manifold, and let $G$ be a compact group acting in a Hamiltonian fashion on $(M,\mathcal{J})$ with generalized moment map $\func{\mu}{M}{\mathfrak{g}^*}$.  Let $a \in \mathfrak{g}^*$, and let $M_a = \bigsqcup (M_{a})_{(H)}$ be the orbit type partition of the GC quotient of $(M,\mathcal{J})$ by $G$ at level $a$.
  Then each component of each $\left( M_{a} \right)_{(H)}$ inherits a natural GC structure from $(M,\mathcal{J})$.

  \item Let $(M,\mathcal{J},H)$ be a \emph{compact} GC manifold, and let $G$ be a compact group acting in a~Hamiltonian fashion on $(M,\mathcal{J},H)$ with generalized moment map $\func{\mu}{M}{\mathfrak{g}^*}$ and moment one-form $\alpha \in \Omega^1(M,\mathfrak{g}^*)$.  Assume that $H+\alpha$ is equivariantly closed.
Let $a \in \mathfrak{g}^*$, and $M_a = \bigsqcup (M_{a})_{(H)}$ be the orbit type partition of the GC quotient of $(M,\mathcal{J},H)$ by~$G$ at level~$a$.
  Then each component of each $\left( M_{a} \right)_{(H)}$ inherits a twisted GC structure from $(M,\mathcal{J},H)$, natural up to $B$-transform.
\end{enumerate}
\end{Theorem}

\begin{proof}
We begin with the twisted case.

First we prove the theorem in the case that $a=0$.

Let $x \in M$ and $K = G_x$.
Note that this implies that $K$ is a closed subgroup of $G$, and is hence compact.
Clearly $K$ acts canonically on $(M,\mathcal{J},H)$, since $G$ does.
By part (c) of Proposition~\ref{prop2}, $M_K^{l_x}$ is open in $M^K$.
It follows that every component of $M^K$ intersecting $M_K^{l_x}$ has the same dimension as $M_K^{l_x}$.
Let $M_x^K$ be the union of components of $M^K$ having nontrivial intersection with $M_K^{l_x}$.
Since each component of $M^K$ is a manifold, it follows that $M_x^K$ is also.
Furthermore, by Proposition~\ref{prop9} we know that each connected component of $M^K$ is a split submanifold of $(M,\mathcal{J})$, and hence also a twisted GC submanifold.
Therefore so is $M_x^K$.

Let $Z_K^x$ be the union of components of $M_K^{l_x}$ that have nontrivial intersection with $\mu^{-1}(0)$.
Since $Z_K^x$ is open in $M_K^{l_x}$, which is open in $M_x^K$, as discussed in Remark~\ref{rem2} we know that $Z_K^x$ is a twisted GC submanifold of $M_x^K$, and hence also of $M$.
Let $\hookfunc{j}{Z_K^x}{M}$ be the inclusion, and denote the ($j^*H$)-twisted GC structure of $Z_K^x$ by $\mathcal{J}'$.

Let $N=N_G(K)$ be the normalizer of $K$ in $G$.  By part (e) of Proposition~\ref{prop2}, we know $M_K^{l_x}$ is $N$-stable.  In fact, so is $Z_K^x$, as we now show.  Note that connected components of manifolds are path-connected.  Therefore, for any $n \in N$, if $y,z \in M_K^{l_x}$ are in the same component, then so are $n \cdot y$ and $n \cdot z$.
Now let $n \in N$ and $y \in Z_K^x$.
By the def\/inition of $Z_K^x$, there exists some $z$ in the same component of $M_K^{l_x}$ as $y$ such that $z \in \mu^{-1}(0)$.
Since $\mu^{-1}(0)$ is $N$-stable, this means $n \cdot z \in \mu^{-1}(0) \cap M_K^{l_x} \subset Z_K^x$, and hence $n \cdot y \in Z_K^x$ as well.

Now we will show that $\mu(Z_K^x)$ and $\alpha( \tang Z_K^x)$ lie in $\operatorname{Ann}_{\mathfrak{g}^*} (\mathfrak{k}) \cap (\mathfrak{g}^*)^K$, where $\operatorname{Ann}_{\mathfrak{g}^*} (\mathfrak{k})$ denotes the annihilator of $\mathfrak{k}$ in $\mathfrak{g}^*$.
Since $M_x^K$ is f\/ixed point-wise by $K$ and $\mu$ and $\alpha$ are equivariant, we know these two sets are contained in $(\mathfrak{g}^*)^K$.
Because $M_x^K$ is closed in $M$, it is compact.  Since $K$ acts trivially on $M_x^K$, it follows from Theorem~\ref{thm6} and Lemmas~\ref{lem6} and \ref{lem1} that $\mathrm{d} \mu^{\xi} = \alpha^{\xi} = 0$ on $\tang M_x^K$, and hence on $\tang Z_K^x$, for all $\xi \in \mathfrak{k}$.
Hence $\mu^{\xi}$ is locally constant on $Z_K^x$ for all $\xi \in \mathfrak{k}$.
Because each component of $Z_K^x$ has nonempty intersection with $\mu^{-1}(0)$, it follows that $\mu^{\xi} = 0$ on $Z_K^x$ for all $\xi \in \mathfrak{k}$, so $\mu(Z_K^x) \subset \operatorname{Ann}_{\mathfrak{g}^*} (\mathfrak{k})$.

Let $L$ denote the quotient Lie group $N/K$, and let $\mathfrak{l}$ denote its Lie algebra.  In Lemma 17 of~\cite{MR1486529}, it is proved that the projection $\mathfrak{g}^* \twoheadrightarrow \mathfrak{n}^*$ dual to the inclusion $\mathfrak{n} \hookrightarrow \mathfrak{g}$ induces an isomorphism
\[
  \operatorname{Ann}_{\mathfrak{g}^*} (\mathfrak{k}) \cap (\mathfrak{g}^*)^K
  \cong
  \operatorname{Ann}_{\mathfrak{n}^*} (\mathfrak{k})
  \cong
  \mathfrak{l}^*.
\]
Let $\func{\mu'}{Z_K^x}{\mathfrak{l}^*}$ and $\func{\alpha'}{\tang Z_K^x}{\mathfrak{l}^*}$ be the compositions of this isomorphism with the restrictions of $\mu$ and $\alpha$, respectively, and note that
\[
  Z_K^x \cap \mu^{-1}(0) = Z_K^x \cap (\mu')^{-1}(0).
\]
Because $Z_K^x$ is f\/ixed point-wise by $K$, the action of $N$ on $Z_K^x$ induces an action of the quotient $L=N/K$ on $Z_K^x$.
We now verify that this action is twisted generalized Hamiltonian with moment map $\mu'$ and moment one-form $\alpha'$.

Since $\mu$, $\alpha$, and the projection $\mathfrak{g}^* \twoheadrightarrow \mathfrak{n}^*$ are $N$-equivariant, and $Z_K^x$ consists of $K$-f\/ixed points, we know that $\mu'$ and $\alpha'$ are $L$-equivariant.  Now we check that $\mu'$ and $\alpha'$ satisfy the generalized moment map conditions for the $L$-action on $Z_K^x$.
Because $K$ f\/ixes the points of $Z_K^x$, the inf\/initesimal action of $\mathfrak{k}$ on $Z_K^x$ is zero, so for all $\xi \in \mathfrak{n}$ we have $[\xi]_{Z_K^x} = \xi_{Z_K^x}$, where $[\xi]$ denotes the image of $\xi$ under the quotient projection $\mathfrak{n} \twoheadrightarrow \mathfrak{n}/\mathfrak{k} \cong \mathfrak{l}$.
As noted above, $\mu^{\eta} = 0$ and $\alpha^{\eta}=0$ for all $\eta \in \mathfrak{k}$, so $(\mu')^{[\xi]} = \mu^{\xi}$ and $(\alpha')^{[\xi]} = \alpha^{\xi}$ for all $\xi \in \mathfrak{n}$.
By Theorem~\ref{thm6} and Lemma~\ref{lem6}, the compositions
  \[
    \xymatrix{
      Z_K^x \ar[r]^{\quad \mu} & \mathfrak{g}^* \ar@{>>}[r] & \mathfrak{n}^*
    }
    \qquad
    \mbox{and}
    \qquad
    \xymatrix{
      \tang Z_K^x \ar[r]^{\quad \alpha}
      & \mathfrak{g}^* \ar@{>>}[r] & \mathfrak{n}^*
    }
  \]
are a generalized moment map and moment one-form for the $N$-action on the ($j^* H$)-twisted GC manifold $Z_K^x$, respectively, so we conclude that
\[
  [\xi]_{Z_K^x} = \xi_{Z_K^x}
  = -\mathcal{J}' \big( \mathrm{d} \mu^{\xi} \big) - \alpha^{\xi}
  = -\mathcal{J}' \big( \mathrm{d} (\mu')^{[\xi]} \big) - (\alpha')^{\xi}
\]
and
\[
  \iota_{[\xi]_{Z_K^x}} (j^* H)
  = \iota_{\xi_{Z_K^x}} (j^* H)
  = \mathrm{d} \alpha^{\xi}
  = \mathrm{d} (\alpha')^{[\xi]}
\]
for all $[\xi] \in \mathfrak{l}$.

By part (e) of Proposition~\ref{prop2}, we know that $N/K$ acts freely on $Z_K^x$, and hence also on~$(\mu')^{-1}(0)$.
Since $H+\alpha$ is $G$-equivariantly closed, we know $H$ is closed.
Using this fact, our computations from the previous paragraph, and part (b) of Remark~\ref{rem1}, we compute
\begin{align*}
  \mathrm{d}_L (j^* H + \alpha')([\xi]) &=
  \mathrm{d} (j^* H) - \iota_{[\xi]_{Z_K^x}} (j^* H)
  + \mathrm{d} (\alpha')^{[\xi]} - \iota_{[\xi]_{Z_K^x}} (\alpha')^{[\xi]} \\
  &= j^*(\mathrm{d} H) - \iota_{\xi_{Z_K^x}} (j^* H)
  + \mathrm{d} \alpha^{\xi} - \iota_{\xi_{Z_K^x}} (\alpha^{\xi}) \\
  &= 0 - \mathrm{d} \alpha^{\xi} + \mathrm{d} \alpha^{\xi} - 0  = 0
\end{align*}
for all $[\xi] \in (\mathfrak{n}/\mathfrak{k})^* \cong \mathfrak{l}^*$.
Hence $j^* H + \alpha'$ is $L$-equivariantly closed.
Therefore we can apply Lin--Tolman's twisted generalized reduction, Theorem~\ref{thm5} above, and obtain a GC structure on the quotient space
\[
  (\mu')^{-1}(0) \big/ (N/K) \cong (\mu')^{-1}(0) \big/ N
  \cong \big(Z_K^x \cap \mu^{-1}(0) \big) \big/ N.
\]
Recall that this structure is only natural up to $B$-transform.  It follows that each component of $\big(Z_K^x \cap \mu^{-1}(0) \big) \big/ N$ is a twisted GC manifold.

By varying the point $x \in M_K$, and thus varying the manifold $Z_K^x$, we can conclude that every component $\big(M_K \cap \mu^{-1}(0) \big) \big/ N$ is a twisted GC manifold.

By parts (d) and (f) of Proposition~\ref{prop2}, we know that $G \cdot M_K = M_{(K)}$ and that the inclusion $M_K \hookrightarrow M_{(K)}$ induces a homeomorphism $M_K/N \approx M_{(K)}/G$.
Together with the fact that $\mu^{-1}(0)$ is $G$-stable, this f\/irst fact implies that $G \cdot \left( M_K \cap \mu^{-1}(0) \right) = M_{(K)} \cap \mu^{-1}(0)$.
Together with the second fact, this implies that
\[
  \left( M_0 \right)_{(K)}
  :=
  \big( M_{(K)} \cap \mu^{-1}(0) \big) \big/ G
  \approx
  \big( M_K \cap \mu^{-1}(0) \big) \big/ N,
\]
and so each component of $\left( M_0 \right)_{(K)}$ inherits a twisted GC structure, natural up to $B$-transform.

The general case, where the reduction is taken at an arbitrary level $a \in \mathfrak{g}^*$ now follows from the shifting trick, as explained following Example~\ref{example5} above.

Now we consider the untwisted case.
Since an untwisted Hamiltonian GC manifold is simply a twisted Hamiltonian GC manifold with $H=0$ and $\alpha = 0$, the only real dif\/ference between parts (a) and (b) of this theorem is that in part (a) we do not assume that $M$ is compact.
Note that the only time above where we used the fact that $M$ is compact was when showing that $\mu^{\xi}$ and $\alpha^{\xi}$ both vanish on $Z_K^x$ for all $\xi \in \mathfrak{k}$, and hence that $\mu(Z_K^x)$ and $\alpha( \tang Z_K^x)$ lie in $\operatorname{Ann}_{\mathfrak{g}^*} (\mathfrak{k})$.
For this non-compact case, note that since $Z_K^x$ contains only $K$-f\/ixed points, we have $\xi_{Z_K^x} = 0$ for all $\xi \in \mathfrak{k}$, so
\[
  \mathrm{d} \mu^{\xi} = \mathcal{J}'(\xi_{Z_K^x})
  = \mathcal{J}'(0) = 0
\]
and hence $\mu^{\xi}$ is locally constant on $Z_K^x$ for all $\xi \in \mathfrak{k}$.
Because each component of $Z_K^x$ has nonempty intersection with $\mu^{-1}(0)$, it follows that $\mu^{\xi} = 0$ for all $\xi \in \mathfrak{k}$, so $\mu(Z_K^x) \subset \operatorname{Ann}_{\mathfrak{g}^*} (\mathfrak{k})$.
This completes the proof of (a).
\end{proof}

\begin{Corollary}[Singular generalized K\"ahler reduction] \label{MainThm2}
The results of Theorem~{\rm \ref{MainThm}} hold if all GC and twisted GC structures are replaced by GK and twisted GK structures, respectively.
\end{Corollary}

\begin{proof}
Suppose $(M,\mathcal{J}_1,\mathcal{J}_2)$ is a GK manifold, twisted or untwisted.  Because a generalized Hamiltonian action on the GK manifold $(M,\mathcal{J}_1,\mathcal{J}_2)$ is simply a generalized Hamiltonian action on the GC manifold $(M,\mathcal{J}_1)$ which also preserves the structure $\mathcal{J}_2$, it is easy to check that the proof of Theorem~\ref{MainThm} holds in precisely the same way for our present situation.  We will simply note that, for any Lie subgroup $K$ of $G$, because both $\mathcal{J}_1$ and $\mathcal{J}_2$ are preserved by $K$, by Proposition~\ref{prop9} we know that each component of $M^K$ is a split submanifold of $M$ with respect to both GC structures, so it is a GK manifold.  Everything else is entirely straightforward to check.
\end{proof}

\subsection*{Acknowledgements}
The author would like to thank Reyer Sjamaar for his help in understanding singular reduction in the symplectic case, Yi Lin for several extremely helpful conversations, Tomoo Matsumara for introducing him to generalized complex geometry, the referees for many useful suggestions, and his family and friends for their unwavering support.

\pdfbookmark[1]{References}{ref}
\LastPageEnding

\end{document}